# Diffusion Based Modeling of Human Brain Response to External Stimuli

Hamidreza Namazi[1], Vladimir V. Kulish[2]

[1,2] School of Mechanical and Aerospace Engineering, Nanyang Technological University, Singapore

[1] Email: m080012@e.ntu.edu.sg   Tel: + 6583740288   (Corresponding Author)

**Abstract**

Human brain response is the overall ability of the brain in analyzing internal and external stimuli in the form of transferred energy to the mind/brain phase-space and thus, making the proper decisions. During the last decade scientists discovered about this phenomenon and proposed some models based on computational, biological, or neuropsychological methods. Despite some advances in studies related to this area of the brain research there was less effort which have been done on the mathematical modeling of the human brain response to external stimuli. This research is devoted to the modeling of human EEG signal, as an alert state of overall human brain activity monitoring, due to receiving external stimuli, based on fractional diffusion equation. The results of this modeling show very good agreement with the real human EEG signal and thus, this model can be used as a strong representative of the human brain activity.

## 1. Introduction

Brain as the most complex organ in the human body controls all bodies' actions/reactions by receiving different stimuli through the nervous system. Any stimulus stronger than the threshold stimulus is translated by the number of sensory neurons generating information about the stimulus and the frequency of the action potentials. After the action potential has been generated, it travels through the neural network to the brain. In various sections of the network and the brain, integration of the signals takes place. Different areas of the brain respond depending on the kind and location of stimuli. The brain sends out signals which generate the response mechanism.

During many years, numerous studies related to brain response to external stimuli have been reported by scientists. Some researchers studied the brain response to different kinds of stimuli without proposing any model. In case of visual stimuli, we can mention the work done by Tomokazu Urakawa et al. in analyzing the effect of the visual stimulus size on the human brain response using Magnetoencephalography (MEG) [1], see also [2-3]. Other group of researchers investigated the effect of auditory stimuli on the brain response. For instance, Udo Will and Eric Berg studied and compared the brain responses to periodic stimulations, silence, and random noise using Electroencephalography (EEG) [4], see also [5-6]. Olfactory stimuli also were the main focus of some researchers. Kouichi Sutani et al. investigated the brain response to pleasant and unpleasant olfactory stimuli using MEG signals. They found out that the MEG signals have recorded from frontal/prefrontal cortical areas of the brain has some differences in case pleasant vs. unpleasant stimuli [7], see also [8-9]. Different works have been done on the brain response to other kinds of stimuli such as emotional stimuli [10-11] and pain stimuli [12-13] also reported in the literature.

On the other hand, some scientist proposed some models of the human brain activity. On a microscopic level, the work done by Freeman in modeling the EEG arising from the olfactory bulb of animals during the perception of odors is note worth. He developed a set of non-linear equations for this response which generates EEG like pattern [14-15]. When the microscopic models are extended to a macro level, then different methods are employed. Many of these models assume the cortical region to be a continuum. Liley et al. developed a set of non-linear continuum field equations which described the macroscopic dynamics of neural activity in the cortical region [16]. These equations were used by Steyn-Ross et al. who introduced noise terms into them to give a set of stochastic partial differential equations (SPDEs). They also converted equations governed by Liley et al. into linearized ODEs. This model could predict the substantial increase in low frequency power at the critical points of induction and emergence. They later used this model to study the electrical activity of an anaesthetized cortex [17-20]. Kramer et al. started with the equations given by Steyn-Ross and co-workers and neglected the spatial variation and the stochastic input. They believed this gave rise to a set of ordinary differential equations (ODEs) for modeling the cortical activity. They showed that the results obtained from the SPDE model agree with clinical data in an approximate way [21], but they also stated that the spatial sampling of the cortex was poor because of inherent shortcomings in the equipment used. Kulish and Chan have suggested a novel method for modeling the brain response by the use of fundamental laws of nature like energy conservation and the least action principle. The model equation obtained has been solved and the results show a good agreement with real EEGs [22].

Robertson et al. claimed that cognitive neuropsychology is more than single-case studies where it can be studied by group design using various methods [23].

Despite rapid advances in studies related to human brain response, there has been less progress in mathematical modeling of human brain activity due to external stimuli. Yet, it seems that the contemporary level of developments in physics and mathematics makes establishing quantitative correlations between external physical stimuli and the brain responses to those stimuli possible.

This paper attempts to introduce a new mathematical model of the human brain response to external stimuli based on fractional diffusion equation.

In this paper first, we talk about the diffusion of energy in human life and then we give a brief description of the macroscopic level of brain organization and EEG as an alert state of human behavior monitoring. Then, by introducing Fractional diffusion equation and the EEG as a multi fractal time series we model the EEG signal using fractional diffusion equation. This model is then solved and by explaining the spectra of fractal dimensions we confirm the multi fractal nature of governed numerical solution. In an application sample, we validate the predicted signal against the real reference EEG signal in case of spectral properties, amplitude and frequency. Some concluding remarks are provided at the end.

**2. Human life and diffusion of energy**

Human life is the total of sensations. Mind is an ability of the Central Nervous System (CNS) that allows human to transform any information about actions or events (external stimuli) into sensations through thoughts. Striving to novel sensations, human strives to activities, the total of which is the maximal for a given individual in the course of his life (Figure. 1).

Evolution of organic & inorganic matter of the earth is anti-entropic; that is energy is stored in structures of increasing complexity. Human is the medium of transforming this energy that otherwise could be stored without any use. There is no limit to the transforming activity of man in principle. As time goes on spatial scale of transformation increases.

Mind is an expression of higher energetic of the CNS in comparison with lower life-forms. Human get an access to energy stored in the planet, practically an unlimited access, much more than is needed for physiological needs, which mean the diffusion of energy from physical phase-space to mind phase-space.

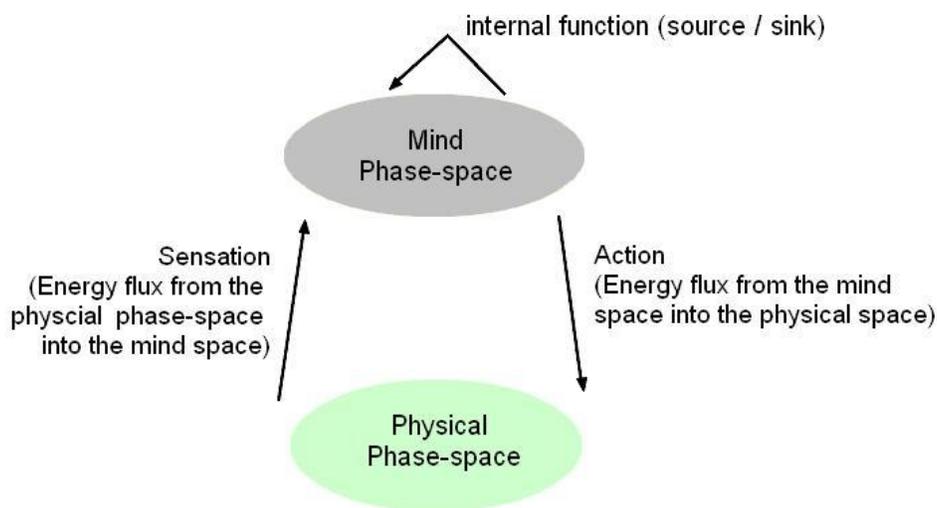

Figure. 1. Conservation of energy in human life

Man is never in equilibrium with the environment. The expanding body of knowledge requires more actions; in order to obtain the same quantity and quality of sensations as it was with a lesser body of knowledge. Actions viewed as environmental changes, require a larger and faster rate of energy transformations. Human actions can be viewed as the diffusion of energy from mind phase-space to physical phase-space. All actions and particularities of man (psyche) and mankind (collective psyche) can be explained by and are determined from the above-mentioned.

On the other hand, the human energy increases/decreases because of some sources/sinks of energy which are not related to his external physical space. These sources/sinks of energy are related to different organisms inside the human body which regulates their tasks by sending/receiving different information (energy) to/from brain.

Considering the above definitions, we can write:

$$\frac{\partial F}{\partial t} = D\nabla^2 F \tag{1}$$

which is the well-known diffusion equation. The coefficient $D$ is the diffusivity of the medium and is closely related to the variance of the random walk performed by the carriers. In fact, this equation arises in many descriptions of biological and physical phenomena, including Brownian motion [24], gradient driven chemical diffusion (with Fick's law), and heat transfer (heat diffusion with Fourier's Law).

It is of great importance to emphasize that the kernel of the diffusion equation of the well-known generalized Gaussian distribution, i.e.,

$$F(r,t) = \frac{1}{(4\pi Dt)^{d/2}} \exp\left(-\frac{r^2}{4Dt}\right) \tag{2}$$

Where $r$ is the radial distance from the origin and $d$ is the Euclidean dimension of the medium in which the process takes place [25].

An alert state of the overall brain activity monitoring device is maintained through monitoring brain signals using EEG patterns. The electric activity of the cerebral cortex highly influences the EEG patterns.

**3. Macroscopic Level of Brain organization**

In order to study the human behavior and in the casual manner neural activity, one can consider the different level of brain organization at many scales in time and space from the single neuron (microscopic level) to the whole brain organization (macroscopic level).

In this research, we focus on the macroscopic level of brain organization. Macroscopic level of brain organization refers to the level of neural assemblies' population in which each neural assembly interacts with other neural assemblies in close and distant cortical areas, exhibits spatial-temporal behaviour, and paint the human behaviour [26]. The functional behaviour of the brain is encoded in these spatial-temporal structures and can be extracted from the macroscopic quantities dynamics observed by EEG signal mostly [27].

Freeman introduced the concept of "wave packet" as a bridge between the microscopic behavior of classical cells and the macroscopic properties of large populations of neurons that have measurable properties such as phase and amplitude [28]. In fact, these wave packets result from synchronization of thousands of neural assemblies in cortical areas of the brain which is reflected by global EEG signals which are the composition of different frequencies (oscillatory activities-Alpha, Beta, …) and are structured coordinatively (spatially-temporally) (Figure. 2). Thus, the global EEG signal has wave-mechanical nature which is proportional to the number of neurons participate in synchronization and is very different from the time-averaged activity of a single neuron.

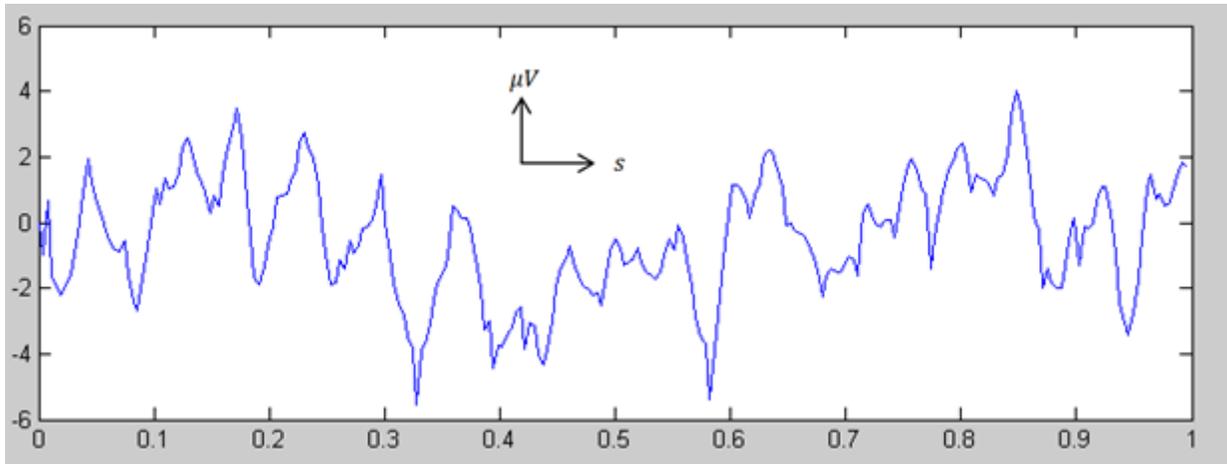

Figure. 2. Global EEG recording - the data are recorded from a normal adult during rest)

As it can be seen in Figure. 2, the EEG signal can be considered as a multi fractal time series which is represented by a series of jumps, randomly directed and of varying length (Brownian motion).

### 4. Fractional diffusion equation

At first glance, the idea to generalize the diffusion process for the case of a fractal space comes when looking at the diffusion equation kernel, equation (2). One is tempted simply to substitute the dimension $d$ by the fractal dimension $\aleph$, and to transform the kernel to a probability distribution which is invariant under addition for $\mathcal{M}$, i.e.,, for which $\sigma^{\aleph} = \sigma_1^{\aleph} + \sigma_2^{\aleph}$ holds (where $\sigma$ is analog to the variance in the Gaussian distribution [29], to obtain the generalized diffusivity.

A very simple relationship exists between these two fractal measures [30]:

$$\langle \aleph \rangle = 1 + d - 2H \tag{3}$$

Where $d$ is the dimension of the Euclidean space in which the fractal process takes place, and $H$ is Hurst exponent has a value within the range $0 < H < 1$ that brings the predictability of signal into account. In fact, the Hurst exponent can be viewed as the probability of the diffusion process being persistent in a certain given direction.

Expressing the Hurst exponent through the representative dimension as:

$$2H = 1 + d - \langle \aleph \rangle \tag{4}$$

One can now use equation (4) to determine the equivalent Hurst exponent when knowing the representative fractal dimension of the domain.

Then, for the transport in a fractal media, the case $H \neq 1/2$ corresponds to an 'anomalous' (fractional) diffusion. Indeed, the diffusion coefficient now may be generalized as:

$$D_\aleph = \frac{1}{2}\frac{d}{dt}\langle B_H(t)^2 \rangle \propto t^{2H-1} \tag{5}$$

Where $B_H(t)$ the fractional Brownian is function of variable $t$ and $\langle B_H(t) \rangle$ denotes the mean of $B_H(t)$ over many samples [25,31]. . Notice that the equality $D = \frac{1}{2}\frac{d}{dt}\langle B_H(t)^2 \rangle$ follows from the relationship between the variance and $D$.

As for the diffusion equation itself, it can be written now in its fractional form as:

$$\frac{\partial^{2H} F}{\partial t^{2H}} = C^{2(2H-1)} D^{2(1-H)} \frac{\partial^2 F}{\partial \eta^2} \tag{6}$$

where $C$ is the speed of propagation. The Hurst exponent, $H$, can be found from equation (4), If the representative fractal dimension of the domain, $\langle \aleph \rangle$, is known. Note that the case $H = 1/2$, which corresponds to a non-fractal diffusion process, leads to the well-known classical equation. Observe also that, if $H = 0$, equation (4) degenerates into the Poisson equation, that is, there is no preferred direction of random walk in this case; while the case $H = 1$ leads to the wave equation.

Equation (6) comes from the substitution of $D_\aleph(t) = CDt^{2H-1}$ into the diffusion equation and then matching the units of both sides of expression, in order to satisfy the dimensional requirement. Taking all these conjectures into account, equation (6) is the only possible form of the equation that is a generalized form for describing the diffusion process when the time becomes fractal.

It is now necessary to make one very important remarks. It is possible to consider a Brownian motion type process as a process which takes place in an Euclidean space [see the right hand side of equation (6)], considering the temporal dimension, $t$, of the process as fractal time which, for the same diffusion coefficient, $D$, either slows down or speeds up the process in question, depending on the Hurst exponent, $H$. This can be described with the generalized diffusion equation (6) in which the time coordinate appears as a fractal quantity. The generalized diffusion equation is a fractional PDE of order $2H$ with respect to time.

## 5. Fractional diffusion model of EEG signal

In this section, we aim to show that EEG as a fractional time series (but not stochastic), which represents a transient record of a random walk process (diffusion process in human life), can be modeled by solutions of the fractional diffusion equation:

$$\frac{\partial^{2H} V}{\partial t^{2H}} = C^{2(2H-1)} D^{2(1-H)} \frac{\partial^2 V}{\partial \eta^2} \tag{7}$$

In order to model EEG signal, here we consider the term $V$ is the voltage fluctuations resulting from ionic current flows within the neurons of the brain which is reflected on the vertical axis in the EEG signal record. The impulse propagation speed is a finite quantity is represented by $C$. The term $D$, diffusion coefficient, as the property of neural tissue, is related to neuron's resistance to the electrical impulse as it travels over the nerve.

One of the interesting points about this model is that it accounts for a finite time lag (reaction time) between any given disturbance (stimulus) and the brain response to it (human action/reaction) based on the assumption that no instantaneous propagation of information is possible within the brain. This effect is considered during the derivation of this model, but it is substituted by $C, D$ according to below equation:

$$\tau = \frac{D}{C^2} \tag{8}$$

In fact, here we consider human behavior, that is reflected in his EEG signal, is the result of energy conservation in human life, and by introducing the diffusion coefficient and its relationship with speed of propagation and time delay, we calculate the value of the diffusion coefficient. This diffusion coefficient stands for the diffusion of energy to human brain in the form of propagated signal from receptor to human brain.

In order to solve the above fractional diffusion model the method proposed by Oldham and Spanier is employed here [32].

Upon introducing the excess value $\hat{V}(\eta,t) = V(\eta,t) - V_0$, so the initial condition imposed on $\hat{V}$ is $\hat{V}(\eta,0) = 0$.

In order to solve the equation (7) we apply the Laplace transform with respect to time, $t$, then we have:

$$C^{2(2H-1)} D^{2(1-H)} \frac{\partial^2 Y}{\partial \eta^2} - s^{2H} Y = 0 \tag{9}$$

Where $s$ is the Laplace transform variable, and $Y$ denotes the Laplace transform of the excess value $\hat{V}$.

The equation (9) is a second-order Ordinary Differential Equation (ODE), where the general solution is calculated as:

$$Y(\eta;s) = A_1(s) e^{-\eta s^H / [C^{(2H-1)} D^{(1-H)}]} + A_2(s) e^{\eta s^H / [C^{(2H-1)} D^{(1-H)}]} \tag{10}$$

Where $A_1(s)$ and $A_2(s)$ are two arbitrary functions. However, since the solution is to be bounded for all $\eta$, the second arbitrary function, $A_2(s)$, must be identically zero, so the solution (10) is changed to:

$$Y(\eta;s) = A_1(s) e^{-\eta s^H / [C^{(2H-1)} D^{(1-H)}]} \tag{11}$$

Upon differentiating (11) with respect to $\eta$:

$$\frac{dY(\eta;s)}{d\eta} = -A_1(s) s^H e^{-\eta s^H / [C^{(2H-1)} D^{(1-H)}]} / [C^{(2H-1)} D^{(1-H)}] \tag{12}$$

After comparing (11) and (12), $A_1(s)$ can be eliminated, then it can be written as:

$$Y(\eta;s) = -s^{-H} C^{(2H-1)} D^{(1-H)} \frac{dY(\eta;s)}{d\eta} \tag{13}$$

By taking the inverse Laplace transform of equation (13) and restoring the original variables then we have:

$$V(\eta,t) = V_0 - C^{(2H-1)} D^{(1-H)} \frac{\partial^{-H}}{\partial t^{-H}} \left(\frac{\partial V}{\partial \eta}\right) \tag{14}$$

which is written in terms of fractional derivative of order $-H$ with respect to $t$.

Using the definition of fractional derivative [32], namely:

$$\frac{\partial^\alpha f}{\partial t^\alpha} = \frac{1}{\Gamma(-\alpha)} \int_0^t \frac{f(\xi)\,d\xi}{(t-\xi)^{\alpha+1}}, \qquad Re(\alpha) < 0 \tag{15}$$

Where $\Gamma(\alpha)$ is the Gamma function, and noticing the Fick's Law:

$$-\frac{\partial V}{\partial \eta} = \frac{\varphi}{D} \tag{16}$$

Where $\varphi$ represents the energy flux, equation (14) can be written as:

$$V(\eta,t) = V_0 + C^{(2H-1)} D^{(-H)} \frac{1}{\Gamma(H)} \int_0^t \frac{\varphi(\eta,\xi)\,d\xi}{(t-\xi)^{1-H}} \tag{17}$$

Since the function $\varphi(\eta,t)$ represents the energy flux through the system boundaries, it can be equated with the external influence acting on the system (external stimulus).

Equation (17) provides the relationship between the non-equilibrium value, $V(\eta,t)$ which stands for the human brain response to any given disturbance (the external influences acting on the human), $\varphi(\eta,t)$. This equation is valid for every location within the domain (including the boundary) at every moment of time.

Since $H$ is a non-negative parameter, it follows from equation (17) that the value of $V$, one the average, increases with the time according to the power law as $t^H$, provided, of course that the fluctuations are small in comparison with the averaged influence. Note the, for $H = 1/2$, $\Gamma(1/2) = \pi^{1/2}$ and the equation (17) yields the well-known diffusion (random walk) growth given by $t^{1/2}$ (see the solution obtained for this case in [32]).

An external influence can be modeled by a Gaussian pulse that is:

$$\varphi(x,t) = \varphi_0(x,t) \exp\left[-\frac{(t-t^*)^2}{\sigma^2}\right] \tag{18}$$

Where $t^*$ denotes the moment of time, at which the Gaussian pulse (external stimulus) reaches its maximal value $\varphi_0(x,t)$, whereas $\sigma$ is the standard deviation of the Gaussian pulse.

In fact, for different type of stimulus, depends on the size and duration, we can have different values of parameters in equation (18).

In the case of many concurrently external stimuli, we will have a series of Gaussian pulses. It is noteworthy to mention that by describing the external influences by a series of Gaussian pulses we consider the same probability of occurrence for each external stimulus,

One very important remark regarding the solution given by equation (17) is to be made at this point. The solution is written in the form of the Volterra integral equation, whose kernel is a power law. It is a well-known that solutions in the form of power law are stable over wide ranges of scale [33]. Hence, the solution for the value evolution on the system given by equation (17) is in general stable. Yet, since many small disturbances may act upon the system even during quite short spans of time, behavior of such a system is never fully stable: the system evolves in a non-equilibrium state.

Another implication from the solution being written in the form of power law is that this solution is self-affine at different scales, i.e. multi fractal. This observation may explain the fact that the pattern of value variations recorded daily is usually similar to the pattern arising in the course of minutes or even hourly recording [29,31], which makes possible use of fractal models for investigating the behavior of EEG signal and accordingly human behavior.

## 6. Spectra of Fractal Dimension

The purpose of this section is to demonstrate that not only do the solutions of fractional diffusion equation resemble fractal time series, but they are indeed multi-fractals, because they do posses non-constant Z-shaped spectra of fractal dimensions.

It has been pointed out [29] that any multi-fractal has to possess a Z-shaped spectrum of fractal dimensions (see also section 14.2 in [25]). A fast algorithm of computing the spectra of a given time series has been developed in [34] and applied to analyze human electroencephalograms.

It has been already pointed out that the concept of fractal dimension is based on the concept of generalized entropy of a probability distribution, introduced by Alfred Renyi [35]. Starting with the moments of order $q$ of the probability $w_i$, Renyi obtained the following expression for entropy:

$$E_q = \frac{1}{1-q} \log \sum_{i=1}^{N} w_i^q \qquad (19)$$

Where $q$ is not necessarily an integer and $log$ denotes $log_2$. Note that for $q \to 1$, equation (19) yields the well-known entropy of a discrete probability distribution [36]:

$$E_1 = -\sum_{i=1}^{N} w_i \log w_i \qquad (20)$$

The probability distribution of a given time series can be recovered by the following procedure. The total range of the value is divided into $N$ bin such that:

$$N = \frac{V_{max} - V_{min}}{\delta V} \tag{21}$$

Where $V_{max}$ and $V_{min}$ are the maximum and the minimum values of the quantity achieved in the course of recording, respectively; $\delta V$ represents the sensitivity (uncertainty) of the determining $V$. The probability that the value falls into the $i$-th bin of size $\delta V$ is computed as:

$$w_i = \lim_{N \to \infty} \frac{N_i}{N} \tag{22}$$

Where $N_i$ equals the number of items the value falls into the $i$-th bin. On the other hand, in case of a time series, the same probability can be found from the ergodic theorem, that is:

$$w_i = \lim_{T \to \infty} \frac{t_i}{T} \tag{23}$$

Where $t_i$ is the time spent by the value in the $i$-th bin during the total time span of recording, $T$.

Further, the generalized fractal dimensions of a given time series with the known probability distribution are defined as:

$$\aleph_q = \lim_{\delta V \to 0} \frac{1}{q-1} \frac{\log \sum_{i=1}^{N} w_i^q}{\log \delta V} \tag{24}$$

Where the parameter $q$ ranges from $-\infty$ to $+\infty$. Note that for a self-similar (simple) fractal time series with equal probabilities $w_i = 1/N$, equation (24) yields $\aleph_q = \aleph_0$ for all values of $q$ [29]. Note also that for a constant value, all probabilities except one become equal to zero, whereas the remaining probability value equals unity. As a result, for a constant value, $\aleph_q = \aleph_0 = 0$. The fractal dimension

$$\aleph_0 = -\frac{\log N}{\log \delta V} \tag{25}$$

is nothing else than Hausdroff-Besicovitch dimension [31].

It is noteworthy to mention that the correlation dimension, frequently used in the analysis of time series, is the fractal dimension with $q = 2$.

As $q \to 1$, equation (24) yields the so-called information dimension:

$$\aleph_1 = \lim_{\delta V \to 0} \frac{-\sum_{i=1}^{N} w_i \log w_i}{\log(1/\delta V)} \tag{26}$$

where the numerator is Shannon's entropy given by (20).

Note also:

$$\aleph_\infty = \lim_{\delta V \to 0} \frac{\log w_{max}}{\log \delta V} \tag{27}$$

and

$$\aleph_{-\infty} = \lim_{\delta V \to 0} \frac{\log w_{min}}{\log \delta V} \tag{28}$$

Such that $\aleph_{-\infty} \geq \aleph_\infty$. In general, if $a < b$, $\aleph_a \geq \aleph_b$, such that $\aleph_q$ is a monotone non-increasing (Z-shaped) function of $\aleph$ [29]. For a given time series ('signal'), the function $\aleph_q$, corresponding to the probability distribution of this series, is called the *fractal spectrum*. Such a name is well-justified, because the fractal spectrum provides information about both frequencies and amplitudes of the series. Indeed, for two probability distributions, a larger value of a fractal dimension of a given order corresponds to the presence of more pronounced spikes (sharper spikes, less expected values of the signal) than in the series for which the value of the fractal dimension of the same order is less. Furthermore, series with a wider range of fractal dimensions, $\aleph_{-\infty} - \aleph_\infty$, can be termed more fractal than signals whose range of fractal dimensions is narrower, so that series with the zero range are self-similar (simple) fractals. In other words, the range of a fractal spectrum is a value associated with the range of frequencies in the series. Now, if the unexpectedness of an event is defined as the inverse of the probability of this event, then steeper spectra correspond to the series in which unexpected values are more dominant, whereas flatter spectra represent those series in which less unexpectedness occurs.

Figure. 3 presents the fractal spectrum of the numerical solution of equation (17). As it can be seen from Figure. 3, this fractal spectrum is a quite widely spread Z-shaped spectrum that corresponds to a multi-fractal time series [29,31]. Hence, the exact solution of the fractional diffusion equation can be used to model multi-fractal EEG time series.

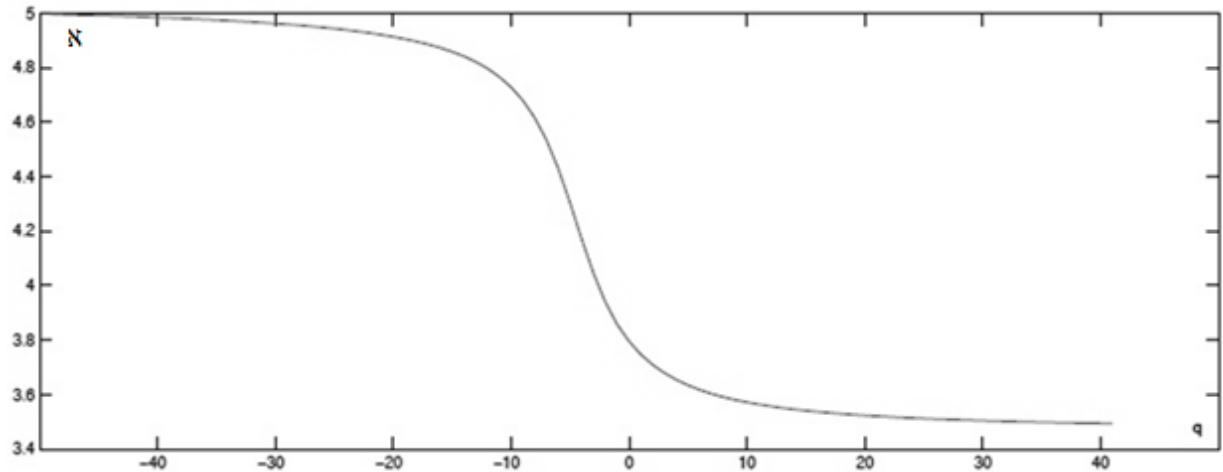

Figure. 3. The fractal spectrum of the numerical solution of the fractional diffusion equation

## 7. Generation of EEG plot by Fractional Diffusion Equation

In order to confirm the usability of the proposed fractional diffusion model, here we apply the model on the real EEG signal and compare the generated signal with the real signal. In fact, by introducing the Hurst exponent in this model we should be able to forecast the EEG signal and accordingly human behavior.

One of the important points about this modeling methodology is that we are not looking for the exact calculated value of signal, like the real reference EEG signal. In fact, we never be able to predict the exact amount of each point in fractal time series like EEG (almost chaotic), but we can know about the strength of prediction by studying the computed values and see how they are similar with the real reference EEG signal in case of amplitude range and frequency (spectral properties analysis). By knowing the amplitude range and frequency of a signal we can study the human behavior. For instance a computed signal with the frequency between 8 to 13 Hz and amplitudes between 5 to 100 microvolt belongs to the brain activity of a normal subject which is relaxing with closed eyes (alpha rhythm).

In order to predict the behavior of EEG fractal time series, we should know about the direction of jumping signal as well as its value. The direction of jump at each point in the EEG signal can be studied by computing the value of a time varying parameter, called Hurst exponent. Hurst exponent is an indicator of the long term memory of the process generating the signal and thus, it is the measure of predictability of signal.

Hurst exponent can have any value between 0 and 1. The value that it gains in each point of signal determines the behavior of next jump in the signal.

In order to calculate the Hurst exponent, we use the algorithm which was formulated by H.E.Hurst [37]. Considering a real EEG signal of a normal human (adult) receiving a single touch stimulus, Figure. 4:

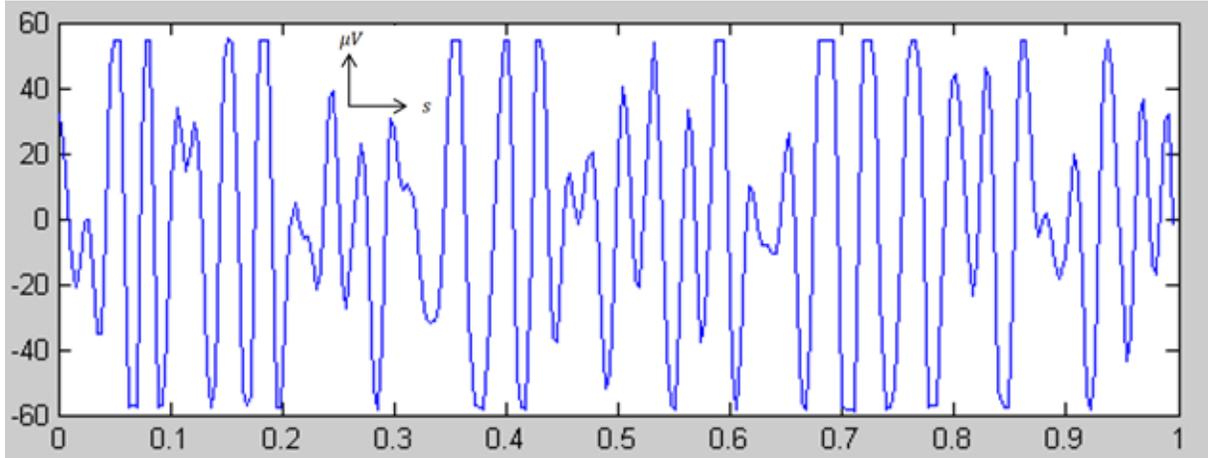

Figure. 4. The reference EEG signal.

With the initial values of $V_0$ =31.99 and $H = 0.79$ at $t = 0$. It is noteworthy to mention that the initial value of $H$ computed from the record of data before t= $0$.

The reference EEG signal has the frequency of 34 Hz and voltage values in the range of:

$$-60\mu V < V < 60\mu V \tag{29}$$

The values to be substituted for the parameters involved in equation (17) for the generation of EEG signal are as given below in Table 1. $C$, the speed of signal propagation is governed from biological studies, and $D$, diffusion coefficient was found using:

$$D = \tau C^2 \tag{30}$$

As we mentioned before $\tau$ accounts for a finite time lag (reaction time) between any given disturbance (stimulus) and the brain response (human action/reaction), based on the assumption that no instantaneous propagation of information is possible within the brain. Reaction time has an average value of 0.215 sec for a normal human.

Table 1. Values of parameters to be substituted in equations (17), and (18)

| Variable | Value | Units |
|---|---|---|
| $c$ | 0.055 | $m/s$ |
| $\tau$ | 0.215 | $s$ |
| $D$ | $6.5 \times 10^{-4}$ | $m^2/s$ |
| $\varphi_0(x,t)$ | 1 | $\frac{V-m}{s}$ |
| $t^*$ | 0.002 | $s$ |
| $C$ | 0.001 | $m/s$ |

Figure. 5 presents the numerical solution of equation (17) for a single touch stimulus on a normal subject which is plotted using MATLAB software. So that the value of generated pulse $V$, is a random variable, uniformly distributed between -60 $\mu V$ and 60 $\mu V$, with the frequency of 34 Hz. It can be easily seen that the solution very closely resembles (in spectral properties) an EEG time series that is shown in Figure 4.

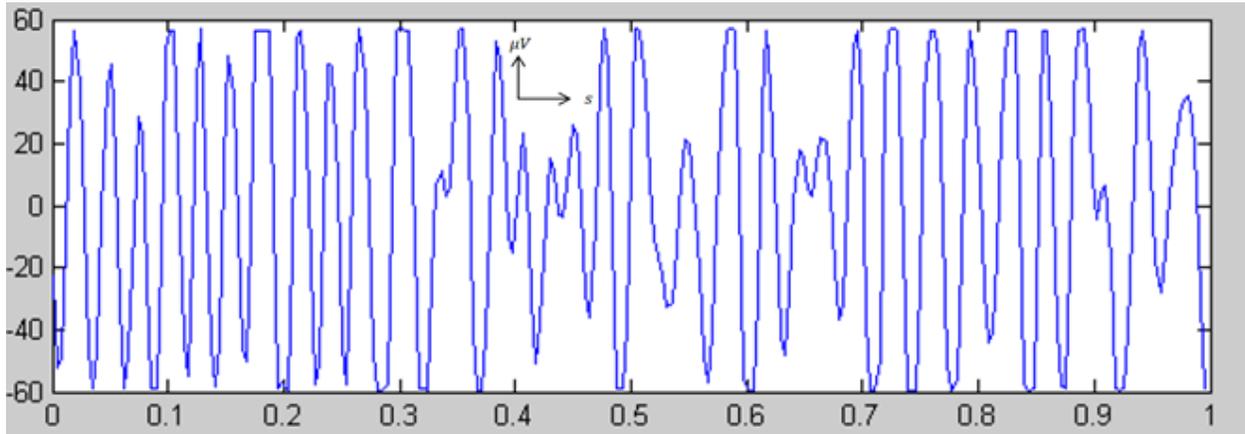

Figure. 5. Numerical solution of equation (17) in the case of an external stimulus.

It can be seen from Figure. 5 that EEG signal can be modeled by the exact solution of fractional partial differential equations and, hence, the behavior of system modeled by means of such equations can, in principle, not only be predicted but also quantifies.

## 8. Conclusion

In this paper we introduced a new mathematical model which represents the human brain response to external stimuli. We developed this model by applying fractional diffusion equation to human EEG signal. The model generates a multi fractal time series which shows a quantitative concurrence in case of spectral properties (amplitude and frequency) with the reference EEG signal. This model shall be further applied in case of different brain activities such as epileptic seizure where the results can be verified against the real EEG signal which means the prediction of the human behavior by forecasting the EEG signal. If so, a seizure warning and the expected time of this epilepsy occurrence can be generated that it means the future monitoring of this disease.